\documentclass[reqno]{amsart}
\usepackage{amsthm, amssymb, amsfonts}
\usepackage{lmodern}
\usepackage{color}
\usepackage{hyperref}
\usepackage[T5]{fontenc}
\usepackage{amscd,amssymb}
\usepackage[v2,cmtip]{xy}
\usepackage{mathrsfs}
\theoremstyle{plain}
\newtheorem{thm}{Theorem}

\newtheorem{corl}[thm]{Corollary}
\theoremstyle{definition}

\theoremstyle{plain}

\theoremstyle{definition}

\allowdisplaybreaks

\begin{document} 
\title[On the hit problem of rank five]
{On the hit problem of rank five}

\author{Nguy\~\ecircumflex n Sum}
\address{Department of Mathematics and Applications, S\`ai G\`on University, 273 An D\uhorn \ohorn ng V\uhorn \ohorn ng, District 5, H\`\ocircumflex\ Ch\'i Minh city, Viet Nam}
 
\email{nguyensum@sgu.edu.vn}

\footnotetext[1]{2000 {\it Mathematics Subject Classification}. Primary 55S10; 55S05.}
\footnotetext[2]{{\it Keywords and phrases:} Steenrod squares, Peterson hit problem, polynomial algebra, modular representation.}

\bigskip
\begin{abstract}
Let $P_q$ be the graded polynomial algebra $\mathbb F_2[x_1,x_2,\ldots ,x_q]$ over the prime field of two elements, $\mathbb F_2$, with the degree of each $x_i$ being 1. We study the {\it hit Peterson problem} of finding a minimal set of generators for $P_q$ as a module over the mod-$2$ Steenrod algebra, $\mathcal{A}.$ In this note, we present the minimal sets of $\mathcal{A}$-generators for $P_5$ in the cases where the degree are sufficiently general.
\end{abstract}
\maketitle

\bigskip
\setcounter{section}{1}
Denote $P_q:= H^*(BE^q) \cong \mathbb F_2[x_1,x_2,\ldots ,x_q],$ a polynomial algebra in  $q$ generators $x_1, x_2, \ldots , x_q$, each of degree 1, where $E^q$ is an elementary abelian 2-group of rank $q$. Here the cohomology is taken with coefficients in the prime field $\mathbb F_2$ of two elements. The algebra $P_q$ is a module over the mod-2 Steenrod algebra, $\mathcal A$.  The action of $\mathcal A$ on $P_q$ is determined by the elementary properties of the Steenrod squares $Sq^i$ and satisfies the Cartan formula (see Steenrod and Epstein~\cite{st}). An element $f$ in $P_q$ is called {\it hit} if it belongs to  $\mathcal{A}^+P_q$, where $\mathcal{A}^+$ is the augmentation ideal of $\mathcal A$.   

We study the {\it Peterson hit problem} of determining a minimal set of generators for the polynomial algebra $P_q$ as a module over the Steenrod algebra. In other words, we want to determine a basis of the $\mathbb F_2$-vector space 
$$QP_q := P_q/\mathcal A^+P_q = \mathbb F_2 \otimes_{\mathcal A} P_q.$$ 

This problem was first studied by Peterson~\cite{pe}, Wood~\cite{wo}, Singer~\cite {si1}, and Priddy~\cite{pr}. Then, this problem  was studied by Carlisle-Wood~\cite{cw}, Crabb-Hubbuck~\cite{ch}, Janfada-Wood~\cite{jw1}, Kameko~\cite{ka}, Mothebe \cite{mo}, Nam~\cite{na}, Repka-Selick~\cite{res}, Silverman~\cite{sl,sl2}, Silverman-Singer~\cite{ss}, Singer~\cite{si2}, Sum-T\'in \cite{su5}, Walker and Wood~\cite{wa3}, Wood~\cite{wo2}, the present author \cite{su1,su2,su4} and others. 

The Peterson hit problem was explicitly determined by 
Peterson~\cite{pe} for $q=1, 2,$ by Kameko~\cite{ka} for $q=3$ and  by us \cite{su2} for $q = 4$. For $q > 4$, this problem is partially determined, All in all, it remains an open question.  The hit problem and its applications to representations of general linear groups have been presented in the monographs of Walker and Wood \cite{wa1,wa2}.

For any nonnegative integer $n$, set $\mu(n) = \min\{k \in \mathbb Z : \alpha (n+k)\leqslant k\}$ where $\alpha (a)$ denotes the number of one in dyadic expansion of a positive integer $a$. 
We denote $(P_q)_n$ and  $(QP_q)_n$ the subspaces of degree $n$ homogeneous polynomials in the spaces $P_q$ and $QP_q$ respectively. The following is Peterson's conjecture, which was established by Wood.

\begin{thm}[See Wood~\cite{wo}]\label{dlmd1} 
	If $\mu(n) > q$, then $(QP_q)_n = 0$.
\end{thm} 

There is a classical operator, known as Kameko's homomorphism
$$\widetilde{Sq}^0_*: (QP_q)_{2d+q} \longrightarrow (QP_q)_d,$$ 
which is induced by an $\mathbb F_2$-linear map $\psi: P_q \to P_q$, given by
$$
\psi(x) = 
\begin{cases}y, &\text{if }x=x_1x_2\ldots x_qy^2,\\  
0, & \text{otherwise,} \end{cases}
$$
for any monomial $x \in P_q$. Note that $\psi$ is not an $\mathcal A$-homomorphism. However, 
$\psi Sq^{2i} = Sq^{i}\psi$ and $\psi Sq^{2i+1} = 0$
for any non-negative integer $i$.

\begin{thm}[See Kameko~\cite{ka}]\label{dlmd2} 
	Let $m$ be a positive integer. If $\mu(2k+q)=q$, then 
	$$(\widetilde{Sq}^0_*)_{(q,k)}:= \widetilde{Sq}^0_*: (QP_k)_{2k+q}\longrightarrow (QP_k)_q$$
	is an isomorphism of the $\mathbb F_2$-vector spaces. 
\end{thm}

Based on Theorems \ref{dlmd1} and \ref{dlmd2}, the hit problem is reduced to the case of degree $n$ of the form
\begin{equation} \label{ct1.1}n =  s(2^d-1) + 2^dm
\end{equation}
with $d, m$ positive integers such that $s-2 \leqslant \mu(m) <s < q$, $\alpha(m + \mu(m)) = \mu(m)$. For $\mu(n)=q-1$, the problem was partially studied by Crabb-Hubbuck~\cite{ch}, Nam~\cite{na}, Repka-Selick~\cite{res}, Walker-Wood \cite{wa3} and the present author ~\cite{su1,su2}. 

\begin{thm}[See \cite{su2}]\label{dlmd3} Let $n =\sum_{1 \leqslant i \leqslant q-1}(2^{d_i}-1)$ 
	with $d_i$ positive integers such that $d_1 > d_2 > \ldots >d_{q-2} \geqslant d_{q-1},$ and let $m = \sum_{1 \leqslant i \leqslant q-2}(2^{d_i-d_{q-1}}-1)$. If $d_{q-1} \geqslant q-1 \geqslant 3$, then
	$$\dim (QP_q)_n = (2^q-1)\dim (QP_{q-1})_m.$$
\end{thm}

For $\mu(n)=q-2$, the hit problem was studied in \cite{su4}. 
Recently, many authors study this problem for the case $q = 5$.

For $q = 5$ the hit problem is reduced to the cases of degrees $n =  s(2^d-1) + 2^dm$ with $1 \leqslant s < 5$ and $\mu(m)<4$. More precisely, we need to consider four families of generic degrees
\begin{align}
&n = 2^{d+1}-1,\ d \geqslant 0,\label{cts1}\\
&n = 2^{d+t} + 2^d - 2\label{cts2}, d > 0,\ t \geqslant 0,\\
&n = 2^{d+t+u} +2^{d+t} + 2^d -3, d,\, u >0,\, t \geqslant 0,\label{cts3}\\
&n = 2^{d+t+u+v} +2^{d+t+u} + 2^{d+t}+ 2^d -4,\ d, u, v > 0,\, t \geqslant 0.\label{cts4}
\end{align}

\bigskip
By studying the structure of $\mathcal A$-generators for $P_5$, we have completely determined the space $(QP_5)_n$ in the cases where the degrees are sufficiently general.

For the case of degree \ref{cts4}, from Theorem \ref{dlmd3} one gets 
\begin{thm}[See \cite{su2}]\label{dlcs40} 
	Let $n = 2^{d+t+u+v} +2^{d+t+u} + 2^{d+t}+ 2^d -4$ with $d,\, t,\, u,\, v$ non-negative integers and $u,\, v > 0$. If $d \geqslant 4$, then
	$$\dim (QP_5)_n = (2^5-1)\dim (QP_{4})_m,$$
	where $m = 2^{t+u+v} +2^{t+u} + 2^{t} -3$.
\end{thm} 

By using Proposition 4.1 in \cite{su2} for $q = 5$, we easily obtain the following.
\begin{thm}\label{dlcs45} Let $n = 2^{d+t+u+v} +2^{d+t+u} + 2^{d+t}+ 2^d -4$ with $d,\, t,\, u,\, v$ positive integers and let $k = 2^{d+t+u+5} +2^{d+t+u} + 2^{d+t}+ 2^d -4$. If $v \geqslant 5$, then $$\dim (QP_5)_n = \dim (QP_5)_{k}.$$
\end{thm}

For $d = 1$, we have

\begin{thm}\label{dlcs41} 
	Let $n_1 = 2^{t+u+v+1} +2^{t+u+1} + 2^{t+1} - 2$ with $t,\, u,\, v$ non-negative integers and $v > 0$. If $t \geqslant 5$, then
	$$\dim (QP_5)_n = 310\dim (QP_{3})_{m},$$
	where $m = 2^{u+v} +2^{u} -2$.
\end{thm}

For $2 \leqslant d \leqslant 3$, we have

\begin{thm}\label{dlcs42} 
	Let $n_2 = 2^{t+u+v+2} +2^{t+u+2} + 2^{t+2}$ with $t,\, u,\, v$ non-negative integers and $v > 0$. If $2 \leqslant d \leqslant 3$ and $t \geqslant 4$, then
	$$\dim (QP_5)_{n_2} = 
	465\dim (QP_{3})_{m},$$
	where $m = 2^{u+v} +2^{u} -2$.
\end{thm} 

For $1 \leqslant d \leqslant 3$ and $t = 0$, we have

\begin{thm}\label{dlcs43} Let $n_0 = 2^{d+u+v} +2^{d+u} + 2^{d+1} - 4$ with $d,\, u,\, v$ positive integers. Then 
	$$\dim (QP_5)_{n_0} = \begin{cases}
	465\dim (QP_{2})_{2^v - 1} + 465\dim (QP_{3})_{2^{v+1}-1}, &\mbox{ if } d = 1,\, u \geqslant 6,\\
	930\dim (QP_{2})_{2^v-1} + 465\dim (QP_{3})_{2^{v+1}-1}, &\mbox{ if } d = 2,\, u \geqslant 5,\\
	1085\dim (QP_{2})_{2^v-1} + 465\dim (QP_{3})_{2^{v+1}-1}, &\mbox{ if } d = 3,\, u \geqslant 5.
	\end{cases}$$
\end{thm}  

For $1 \leqslant d \leqslant 3$ and $0 < t < 5$, by using Theorem \ref{dlcs2} below we easily obtain the following.
\begin{thm}\
	Let $n = 2^{d+t+u+v} +2^{d+t+u} + 2^{d+t}+ 2^d -4$ with $d,\, t,\, u,\, v$ positive integers and let $k = 2^{d+t+v+6} +2^{d+t+6} + 2^{d+t}+ 2^d -4$. If $u \geqslant 6$, then
	$$\dim (QP_5)_n = \dim (QP_{5})_k.$$
\end{thm}

\bigskip
For the case of degree \ref{cts3}, $n = 2^{d+t+u} +2^{d+t} + 2^d -3$, we consider Kameko's homomorphism
$$(\widetilde{Sq}^0_*)_{(5,n)}: (QP_5)_{n}\to (QP_5)_{\frac{n-5}2}.$$
Here $\frac{n-5}2 = 2^{d+t+u-1} +2^{d+t-1} + 2^{d-1} -4$.

Since Kameko's homomorphism is an epimorphism, by combining this and the result for $s = 4$, we get

\begin{thm}[See \cite{su4}] Let $n = 2^{d+t+u} + 2^{d + t} + 2^d -3$ with $d,\, t,\, u$ positive integers. If  $d \geqslant 6$ and $t,\, u \geqslant 4$, then
	$$\dim (QP_5)_n = 4(2^3-1)(2^4-1)(2^5-1) = 13020.$$
\end{thm}

Recently, we have proved the following.

\begin{thm}\label{dlcs3} Let $n = 2^{d+t+u} + 2^{d + t} + 2^d -3$ with $d,\, t,\, u$ non-negative integers and $u > 0$. If  $d \geqslant 5$, then
	$$\dim(\mbox{\sf Ker}(\widetilde{Sq}^0_*)_{(5,n)}) = 155 \dim(QP_3)_m,$$
	where $m = 2^{t+u} + 2^{t} -2$.
\end{thm}

By combining Theorems \ref{dlcs40} and \ref{dlcs3} one gets the following.

\begin{corl}\label{hhq} Let $n = 2^{d+t+u} + 2^{d + t} + 2^d -3$, with $d,s,t$ integers such that $d \geqslant 5,\, t \geqslant 0$ and $u > 0$. The dimension of the $\mathbb F_2$-vector space $(QP_5)_n$ is given by the following table: 
	
	\medskip
	\centerline{\begin{tabular}{c|cccccc}
			$n$ & $u=1$ & $u=2$ & $u=3$ &$u=4$ &$u=5$  & $u\geqslant 6$\cr
			\hline
			$t = 0$ & \textcolor{blue}{$1116$} & \textcolor{blue}{$2790$} & $3813$ & $4960$ & $5735$ & $6045$ \cr
			$t = 1$ & $3410$ & $6231$ & $7285$ & $7719$ & $7595$ & $7595$ \cr
			$t = 2$ & $5766$ & $9207$ & $10726$ & $11160$ & $11160$ & $11160$ \cr
			$t = 3$ & $7254$ & $10695$ & $12090$ & $12555$ & $12555$ & $12555$ \cr
			$t \geqslant 4$ & $7595$ & $11160$ & $12555$ & $13020$ & $13020$ & $13020$ \cr
	\end{tabular}}
\end{corl}
We note that the results for the cases $t = 0,\, u = 1$ and $t = 0,\, u = 2$ had been announced in \cite{ph} and \cite{subs}. 

\medskip
For the case $d = 1$, we have 

\begin{thm}\label{dlcs31} Let $n = 2^{t+u+1} + 2^{t+1} - 1$ with $t,\, u$ positive integers. If  $t \geqslant 5$, then
	$$\dim(\mbox{\sf Ker}(\widetilde{Sq}^0_*)_{(5,n)}) = 465 \dim(QP_2)_{2^u-1}.$$
\end{thm}

For the case $d = 2$, we have 

\begin{thm}\label{dlcs32} Let $n = 2^{t+u+2} + 2^{t+2} + 1$ with $t,\, u$ non-negative integers. If  $t \geqslant 4$, then
	$$\dim(\mbox{\sf Ker}(\widetilde{Sq}^0_*)_{(5,n)}) = 930 \dim(QP_2)_{2^u-1}.$$
\end{thm}

For the case $d = 3,\, 4$, we have 

\begin{thm}\label{dlcs33} Let $n = 2^{d+t+u} + 2^{d+t} + 2^{d} - 3$ with $t,\, u$ positive integers. If $3\leqslant d \leqslant 4$ and $t \geqslant 3$, then
	$$\dim(\mbox{\sf Ker}(\widetilde{Sq}^0_*)_{(5,n)}) = 1085 \dim(QP_2)_{2^u-1}.$$
\end{thm}

\begin{thm}\label{dlcs30} Let $n = 2^{d+t+u} + 2^{d+t} + 2^{d} - 3$ with $t,\, u$ positive integers and let $k = 2^{d+t+5} + 2^{d+t} + 2^{d} - 3$. If $u \geqslant 5$, then
	$$\dim(\mbox{\sf Ker}(\widetilde{Sq}^0_*)_{(5,n)}) = \dim(\mbox{\sf Ker}(\widetilde{Sq}^0_*)_{(5,k)}).$$
\end{thm}

For the case of degree \ref{cts2}, $n = 2^{d+t} +2^{d} - 2$, we have
\begin{thm}\label{dlcs2} Let $n = 2^{d+t} + 2^{d} - 2$ with $d,\, t$ non-negative integers. If  $d \geqslant 6$, then
	\begin{align*}\dim (QP_5)_n &= 155\dim(QP_2)_{2^t-1} + 310\dim(QP_3)_{2^{t+1}-1}.
	\end{align*}
\end{thm}
By using the resuls of Peterson \cite{pe} and Kamekp \cite{ka} we obtain the following.

\begin{corl} Let $n = 2^{d+t} + 2^{d} - 2$ with $d,\, t$ non-negative integers and $d \geqslant 6$. Then, we have
	\begin{align*}\dim (QP_5)_n = \begin{cases} \textcolor{blue}{1085}, &\mbox{if } t = 0,\\
	2480, &\mbox{if } t = 1,\\
	3565, &\mbox{if } t = 2,\\
	4495, &\mbox{if } t = 3,\\
	4805, &\mbox{if } t \geqslant 4.
	\end{cases}
	\end{align*}
\end{corl}
The result for the case $t = 0$ had been announced in \cite{sux}. 

\medskip
For the case $1 \leqslant d \leqslant 5$, we have

\begin{thm}\label{dlcs20} Let $n = 2^{d+t} + 2^{d} - 2$ with $d,\, t$ positive integers and let $k = 2^{d+5} + 2^{d} - 2$. If $t \geqslant 5$, then
	$$\dim (QP_5)_n = \dim (QP_5)_k.$$
\end{thm} 

\medskip
For the case of degree \ref{cts1}, $n = 2^{d+1} -1$, with $d \geqslant 5$, we consider Kameko's homomorphism
$$(\widetilde{Sq}^0_*)_{(5,n)}: (QP_5)_{(2^{d+1}-1)}\longrightarrow (QP_5)_{(2^{d}-3)}.$$
This is an epimorphism and the space $(QP_k)_{(2^{d}-3)}$ has been computed in \cite{ph}. So, we need only to determine $\mbox{\rm Ker}\big((\widetilde{Sq}^0_*)_{(5,2^{d}-3)}\big)$.

\begin{thm}[See \cite{sux}]\label{thm2}
	For any integer $d \geqslant 6$, there exist exactly $496$ classes of degree $2^{d+1} -1$ in $\mbox{\rm Ker}\big((\widetilde{Sq}^0_*)_{(5,n)}\big)$ represented by the admissible monomials in $P_5$. Consequently 
	$\dim \mbox{\rm Ker}\big((\widetilde{Sq}^0_*)_{(5,2^{d}-3)}\big) = 496.$
\end{thm}

Combining Theorem \ref{thm2} and the results in Ph\'uc \cite{ph}, we obtain the following. 

\begin{corl}[See \cite{sux}]\label{hqdl2} For any integer $d \geqslant 6$, we have
	$$\dim(QP_5)_{(2^{d+1}-1)} = \begin{cases}1441, \mbox{ if } d = 6,\\ 1611, \mbox{ if } d = 7,\\ 1612, \mbox{ if } d \geqslant 8. \end{cases} $$
\end{corl}

The results of this paper show that to complete the hit problem of rank five, it suffices to determine $(QP_5)_{n}$ for some special cases of degree $n$, taking the following forms:
\begin{align}
	&n = 2^{d+t} + 2^d - 2,\ d \leqslant 5,\ t \leqslant 5,\\
	&n = 2^{d+t+u} +2^{d+t} + 2^d -3,\ d \leqslant 4,\, t \leqslant 4,\, u \leqslant 5,\\
	&n = 2^{d+t+u+v} +2^{d+t+u} + 2^{d+t}+ 2^d -4,\ d \leqslant 3,\, t \leqslant 4,\, u,\, v \leqslant 6.
\end{align}
The notions related to the hit problem and the notations are used following \cite{su2} and \cite{su4}.

For $J= (j_1, j_2, \ldots, j_u) : 1 \leqslant j_1 <\ldots < j_u \leqslant q$, we define a monomorphism $\theta_J: P_u \to P_q$ of $\mathcal A$-algebras by substituting 
\begin{equation}\label{ctbs}
\theta_J(x_s) = x_{j_s} \ \mbox{ for } \ 1 \leqslant s \leqslant u.
\end{equation} 
It is easy to see that, for any weight vector $\omega$ of degree $n$, 
\[Q\theta_J(P_u^+)(\omega) \cong  QP_u^+(\omega)\mbox{ and } (Q\theta_J(P_u^+))_n \cong (QP_u^+)_n\] 
for $1 \leqslant u \leqslant q$, where $Q\theta_J(P_u^+) = \theta_J(P_u^+)/\mathcal A^+\theta_J(P_u^+)$. So, by a simple computation, we get the following which is presented Walker and Wood~\cite{wa1}.

For a weight vector $\omega$ of degree $n$, we have direct summand decompositions of the $\mathbb F_2$-vector spaces
\begin{align*} QP_q(\omega)  &= \bigoplus_{\mu(n) \leqslant u\leqslant q}\bigoplus_{\ell(J) =u}Q\theta_J(P_u^+)(\omega), 
\end{align*}
where $\ell(J)$ is the length of $J$. Consequently 
\begin{align}
\dim QP_q(\omega) &= \sum_{\mu(n) \leqslant u\leqslant q}{q\choose u}\dim QP_u^+(\omega),\label{ctc1}\\
\dim (QP_q)_n &= \sum_{\mu(n) \leqslant u\leqslant q}{q\choose u}\dim (QP_u^+)_n.\label{ctc2}
\end{align}

We present the set of admissible monomials with degree forms from \ref{cts1} to \ref{cts4}. The notions related to the hit problem and the notations are used following \cite{su2} and \cite{su4}. The set of admissible monomials of degrees \ref{cts1}, \ref{cts3}, \ref{cts4} are explicitly determined in \cite{sux}, \cite{su4}, \cite{su2} respectively. So we need only to present the set of admissible monomials of degree \ref{cts2}. Let $n = 2^{d+t} + 2^d - 2$ with $d \geqslant 6$. We see that for any $x \in B_5((2)|^{d})$ there exists uniquely a pair $\mathbb J^x = (j_1^x,j_2^x)$ such that $1 \leqslant j_1^x<j_2^x \leqslant 5$ and $\nu_{j_1^x}(x)>16$, $\nu_{j_2^x}(x) > 16$. For any $x \in B_5((4)|(3)|^{d-2})$ there exists uniquely a triple $\tilde{\mathbb J}^x = (\tilde j_1^x,\tilde j_2^x,\tilde j_3^x)$ such that $1 \leqslant \tilde j_1^x < \tilde j_2^x < \tilde j_3^x \leqslant 5$ and $\nu_{\tilde j_1^x}(x)>16$, $\nu_{\tilde j_2^x}(x) > 16$, $\nu_{\tilde j_3^x}(x) > 16$. We prove Theorem \ref{dlcs2} by proving the following.

\begin{align*}
B_5(n) &= \left\{x\theta_{\mathbb J^x}\left(z^{2^{d}}\right): x \in B_5((2)|^{d}),\, z \in B_2(2^t-1)\right\}\\
&\quad \bigcup\left\{x\theta_{\tilde {\mathbb J}^x} \left(z^{2^{d-1}}\right): x \in B_5((4)|(3)|^{d-2}),\, z \in B_3(2^{t+1}-1
\right\}.
\end{align*}

The sets $B_5((2)|^d)$ and $B_5((4)|(3)|^{d-2})$ are explicitly determined in \cite{sux}.

Some of the results of this paper were presented at the 2023 Vietnam Mathematical Congress in Da Nang. The detailed proofs of the results will appear in the near future. Some other new results will be updated in the next version.

{}


\begin{thebibliography}{99}

\bibitem{cw}  D. P. Carlisle and R. M. W. Wood, \textit{The boundedness conjecture for the action of the Steenrod algebra on polynomials}, in: N. Ray and G. Walker (ed.), Adams Memorial Symposium on Algebraic Topology 2, (Manchester, 1990), in: London Math. Soc. Lecture Notes Ser., Cambridge Univ. Press, Cambridge, vol. 176,  1992, pp. 203-216,  MR1232207. 

\bibitem{ch}  M. C. Crabb and J. R. Hubbuck, \textit{Representations of the homology of $BV$ and the Steenrod algebra II}, Algebraic Topology: new trend in localization and periodicity, Progr. Math. 136 (1996) 143-154,  MR1397726. 

\bibitem{jw1} A. S. Janfada and R. M. W. Wood, \textit{The hit problem for symmetric polynomials over the Steenrod algebra}, Math. Proc. Cambridge Philos. Soc. 133 (2002) 295-303, MR1912402.

\bibitem{ka}  M. Kameko, \textit{Products of projective spaces as Steenrod modules}, PhD. Thesis, The Johns Hopkins University, ProQuest LLC, Ann Arbor, MI, 1990. 29 pp., MR2638633. 

\bibitem{na} T. N. Nam, \textit{$\mathcal{A}$-g\'en\'erateurs g\'en\'eriquess pour l'alg\`ebre polynomiale}, Adv. Math. 186 (2004) 334-362, MR2073910. 
 
\bibitem{pe}  F. P. Peterson, \textit{Generators of $H^*(\mathbb RP^\infty \times \mathbb RP^\infty)$ as a module over the Steenrod algebra},  Abstracts Amer. Math. Soc. No. {833} (1987) 55-89.

\bibitem{ph} \DJ. V. Ph\'uc, \textit{The "hit" problem of five variables in the generic degree and its application}, Topology Appl. 282 (2020) 107321,  MR4123276, https://doi.org/10.1016/j.topol.2020.107321.

\bibitem{pr}  S. Priddy, \textit{On characterizing summands in the classifying space of a group, I}, Amer. Jour. Math. 112 (1990) 737-748, MR1073007. 

\bibitem{res}  J. Repka and P. Selick, \textit{On the subalgebra of $H_*((\mathbb RP^\infty)^n;\mathbb F_2)$ annihilated by Steenrod operations}, J. Pure Appl. Algebra 127 (1998) 273-288, MR1617199.

\bibitem{sl}   J. H. Silverman, \textit{Hit polynomials and the canonical antiautomorphism of the Steenrod algebra}, Proc. Amer. Math. Soc. 123 (1995) 627-637, MR1254854. 

\bibitem{sl2}   J. H. Silverman, \textit{Hit polynomials and conjugation in the dual Steenrod algebra}, Math. Proc. Cambridge Philos. Soc. 123 (1998) 531-547, MR1607993. 

\bibitem{ss}   J. H. Silverman and W. M. Singer, \textit{On  the  action  of Steenrod  squares  on  polynomial  algebras  II},  J. Pure Appl. Algebra 98 (1995) 95-103, MR1317001.

\bibitem{si1}  W. M. Singer, \textit{The transfer in homological algebra}, Math. Zeit.  202 (1989) 493-523, MR1022818. 

\bibitem{si2}   W. M. Singer, \textit{On the action  of the Steenrod squares on polynomial algebras}, Proc. Amer. Math. Soc.  111 (1991) 577-583, MR1045150. 

\bibitem{st}  N. E. Steenrod and D. B. A. Epstein, \textit{Cohomology operations},   Annals of Mathematics Studies 50,  Princeton University Press, Princeton N.J (1962), MR0145525. 

\bibitem{su1}  N. Sum, \textit{The negative answer to Kameko's conjecture on the hit problem}, Adv. Math. {225} (2010)  2365-2390, MR2680169.

\bibitem{su2} N. Sum, \textit{On the Peterson hit problem}, Adv. Math. 274 (2015) 432-489, MR3318156.

\bibitem{su4} N. Sum, \textit{The squaring operation and the hit problem for the polynomial algebra in a type of generic degree}, J. Algebra 622 (2023) 165-196, MR4547877.

\bibitem{sux} N. Sum, \textit{The admissible monomial bases for the polynomial algebra of five variables in some types of generic degrees}, Topology Appl. 349 (2024), 108909.

\bibitem{subs}  N. Sum, \textit{The hit problem of rank five in a generic degree}, East-West J. Math. 25 (2024), no. 1, 1-60, MR4750740.

\bibitem{su5} N. Sum and N. K. T\'in, \textit{The hit problem for the polynomial algebra in some weight vectors}, Topology Appl. 290 (2021) 107579, MR4199843, https://doi.org/10.1016/j.topol.2020.107579.

\bibitem{sut} N. Sum, P. D Tai, \textit{On a minimal set of generators for the polynomial algebra of five variables in a generic degree}, Asian-Eur. J. Math., 18 (2025), no. 12, Paper No. 2550055, 38 pp. MR4982164

\bibitem{wa3}  G. Walker and R. M. W. Wood, \textit{Flag modules and the hit problem for the Steenrod algebra}, Math. Proc. Cambridge Philos. Soc.  147 (2009) 143-171, MR2507313. 

\bibitem{wa1} G. Walker and R. M. W. Wood, \textit{Polynomials and the mod 2 Steenrod algebra, Vol. 1: The Peterson hit problem}, London Mathematical Society Lecture Note Series 441, Cambridge University Press, 2018,  MR3729477. 

\bibitem{wa2} G. Walker and R. M. W. Wood, \textit{Polynomials and the ${\rm mod}\, 2$ Steenrod algebra. Vol. 2. Representations of ${\rm GL}(n,\mathbb F_2)$}, London Mathematical Society Lecture Note Series, 442. Cambridge University Press, 2018, MR3729478.

\bibitem{wo}  R. M. W. Wood, \textit{Steenrod squares of polynomials and the Peterson conjecture}, Math. Proc. Cambriges Phil. Soc. 105  (1989) 307-309, MR0974986. 

\bibitem{wo2}  R. M. W. Wood, \textit{Problems in the Steenrod algebra}, Bull. London Math.  Soc.  30 (1998) 449-517, MR1643834. 

\end{thebibliography}
\end{document}